\documentclass[10pt,a4paper]{article}

\usepackage{graphicx,amssymb}

\usepackage{amsmath}
\usepackage{amsfonts}
\usepackage{verbatim}
\usepackage{xy}
\xyoption{all}
\CompileMatrices

\renewcommand{\i }{\iota}

\newcommand{\lgs}[2]{\mbox{$\bigl(\frac{#1}{#2}\bigr)$}}
\newcommand{\isom}{\simeq}
\newcommand{\LR}{\longrightarrow}
\newcommand{\gothr}{{\mathfrak r}}

\DeclareMathOperator{\Gal}{Gal}
\DeclareMathOperator{\rk}{rk}

\DeclareMathOperator{\Res}{Res}

\DeclareMathOperator{\NO}{{\cal N}}
\DeclareMathOperator{\NCO}{{\lceil{\cal N}\rceil}}
\DeclareMathOperator{\N}{{\NO\!}}
\DeclareMathOperator{\NC}{{\NCO\!}}

\newcommand{\ov}[1]{\overline{#1}}
\newcommand{\Q}{{\mathbb Q}}
\newcommand{\Z}{{\mathbb Z}}
\newcommand{\F}{{\mathbb F}}

\newcommand{\z}{{\zeta}}

\renewcommand{\th}{\theta}
\newcommand{\om}{\omega}

\renewcommand{\a}{{\mathfrak a}}
\renewcommand{\b}{{\mathfrak b}}
\renewcommand{\c}{{\mathfrak c}}
\newcommand{\f}{{\mathfrak f}}
\newcommand{\al}{\alpha}
\newcommand{\be}{\beta}
\newcommand{\ga}{\gamma}

\newcommand{\p}{{\mathfrak p}}
\newcommand{\q}{{\mathfrak q}}
\newcommand{\gd}{{\mathfrak d}}

\newcommand{\eps}{\varepsilon}

\newcommand{\myproof}[1]{{\it Proof #1. \/}}
\newcommand{\Proof}{{\it Proof. \/}}

\newcommand{\squareforqed}{\hbox{\rlap{$\sqcap$}$\sqcup$}}
\newcommand{\qed}{\ifmmode\squareforqed\else{\unskip\nobreak\hfil
\penalty50\hskip1em\null\nobreak\hfil\squareforqed
\parfillskip=0pt\finalhyphendemerits=0\endgraf}\fi}

\newcommand{\fp}{\qed\removelastskip\vskip\baselineskip\relax}

\renewcommand{\pmod}[1]{\allowbreak\ ({\rm{mod}}\,\,#1)}

\newtheorem{theorem}{Theorem}[section]
\newtheorem{corollary}[theorem]{Corollary}
\newtheorem{proposition}[theorem]{Proposition}
\newtheorem{lemma}[theorem]{Lemma}
\newtheorem{definition}[theorem]{Definition}

\begin{document}
 \pagestyle{plain}

\title{Counting Cubic Extensions with given Quadratic Resolvent}
\author{Henri Cohen and Anna Morra\footnote{The second author was supported by
the European Community under the Marie Curie Research Training Network GTEM 
(MRTN-CT-2006-035495)} ,\\
Universit\'e Bordeaux I, Institut de Math\'ematiques de Bordeaux,\\
351 Cours de la Lib\'eration, 33405 TALENCE Cedex, FRANCE}
\maketitle

\begin{abstract}
Given a number field $k$ and a quadratic extension $K_2$, we give an explicit
asymptotic formula for the number of isomorphism classes of cubic extensions
of $k$ whose Galois closure contains $K_2$ as quadratic subextension, ordered
by the norm of their relative discriminant ideal. The main
tool is Kummer theory. We also study in detail the error term of the
asymptotics and show that it is $O(X^{\alpha})$, for an explicit $\alpha<1$.
\end{abstract}


\section{Introduction and Statement of Results}

\subsection{Introduction}

Let $k$ be a number field, fixed once and for all as our base field, let
$K/k$ be a cubic extension of $k$, and let $N$ be a Galois
closure of $K/k$. When $K/k$ is not cyclic we have
$\Gal(N/k)\isom S_3\isom D_3$, and the field $N$ contains a unique quadratic subextension
$K_2/k$. 

When $K/k$ is cyclic we have $N=K$ and $\Gal(N/k)\isom C_3$. Although this
case has already been treated in \cite{CoDiOl3}, since the methods are
almost identical we include it in the present paper by setting $K_2=k$,
which by abuse of language we will still call a quadratic extension of $k$,
even though $[K_2:k]=1$.

We fix the quadratic extension $K_2/k$, and we call $\mathcal{F}(K_2)$ the set
of cubic extensions $K/k$ (up to $k$-isomorphism) such that the 
quadratic subextension of the Galois closure of $K/k$ is isomorphic to $K_2$.
Our goal is to compute an asymptotic formula for
$$N(K_2/k,X)=|\{K\in\mathcal{F}(K_2),\ \N_{k/\Q}(\gd(K/k))\le X\}|\;,$$
where $\gd(K/k)$ is the relative discriminant ideal of $K/k$ and 
$\N_{k/\Q}$ denotes the absolute norm.

By a well-known theorem (see for example Theorem 9.2.6 of \cite{Coh2}), the
conductor of the cyclic extension $N/K_2$ is of the form
$\f(N/K_2)=\f(K/k)\Z_{K_2}$, where $\f(K/k)$ is an ideal of the base field $k$
(when $K/k$ is noncyclic this is of course not a conductor in the usual
sense). When $k=\Q$ we will write $f(K)$ for the positive integer generating
the ideal $\f(K/\Q)$ of $\Z$.

Since $\gd(K/k)=\gd(K_2/k)\f(K/k)^2$, it is clear that
$$N(K_2/k,X)=M(K_2/k,(X/\N_{k/\Q}(\gd(K_2/k)))^{1/2})\;, \textrm{ where}$$
$$M(K_2/k,X)=|\{K\in\mathcal{F}(K_2),\ \N_{k/\Q}(\f(K/k))\le X\}|\;,$$
so we will in fact only study $M(K_2/k,X)$. When $k=\Q$, we will omit the
letter $k$ from the notation.

Some results of this paper are obtained using tools which are similar
to the ones used (in a slightly different context) in previous papers
of the first author and collaborators (\cite{CoDiOl2}, \cite{CoDiOl3}).
Thus, for brevity we have decided to omit or only sketch some long and 
technical proofs, and we refer to (\cite{CoDiOl2}, \cite{CoDiOl3})
for complete proofs which can be easily adapted to our situation.

On the other hand we would like to emphasize that the Galois structure
and the use of Kummer theory are more complex in our case than the in the 
cyclic case (\cite{CoDiOl3}), so some results require new proofs, which we 
give in detail.

Moreover, unlike (\cite{CoDiOl3}), we give an explicit formula for
the error term, since this kind of technique, although considered 
``standard'', is not easy to find in detail in the literature.

Finally, in some cases it is possible to give simple explicit formulas by
using Scholz's Spiegelungssatz, and this is done in Section \ref{schol}.

\subsection{Statement of Results}

The result in the case of a general base field $k$ is a little complicated
(see Corollary \ref{resgen1}), so we state it here only for $k=\Q$.

\begin{theorem}\label{thresq} As above, let $K_2=\Q(\sqrt{D})$ be an extension
of $\Q$ with $[K_2:\Q]\le2$, denote by $K'_2=\Q(\sqrt{-3D})$ the mirror field
of $K_2$, and set $g(K'_2)=3$ if $K'_2=\Q(\sqrt{-3})$, and $g(K'_2)=1$ 
otherwise. Then:
\begin{enumerate}\item (Pure cubic fields.) We have
$$M(\Q(\sqrt{-3}),X)=C(\Q(\sqrt{-3}))X(\log(X)+D(\Q(\sqrt{-3}))-1)+O(X^{2/3+\varepsilon}),$$
for every $\varepsilon>0$, where
\begin{align*}
C(\Q(\sqrt{-3}))&=\dfrac{7}{30}\prod_{p}\left(1-\dfrac{3}{p^2}+\dfrac{2}{p^3}\right)\\
D(\Q(\sqrt{-3}))&=2\ga-\dfrac{16}{35}\log(3)+6\sum_p\dfrac{\log(p)}{p^2+p-2}\;,\end{align*}
and $\ga$ is Euler's constant.
\item (General case.) For $D\ne-3$, denote by $a_{K'_2}(p)$ the number of 
copies of $\Q_p$ occurring in $K'_2\otimes\Q_p$ ($a_{K'_2}(p)=0$ or $2$
according to whether the number of prime ideals above $p$ in $K'_2$ is equal
to $1$ or $2$). Then
$M(\Q(\sqrt{D}),X)=C(\Q(\sqrt{D}))X+O(X^{2/3+\varepsilon})$, where
$$C(\Q(\sqrt{D}))=g(K'_2)\dfrac{c_3(K'_2)}{3^{3+r_2(K'_2)}}\prod_{p\ne 3}\left(1+\dfrac{a_{K'_2}(p)}{p}\right)\left(1-\dfrac{1}{p}\right)\;,$$
and
$$c_3(K'_2)=\begin{cases}
11&\text{\quad if $3\Z_{K'_2}=\p_1^2$\;,}\\
15&\text{\quad if $3\Z_{K'_2}=\p_1$\;,}\\
21&\text{\quad if $3\Z_{K'_2}=\p_1\p_2$\;.}
\end{cases}$$
\end{enumerate}\end{theorem}

The result of (2) for $D=1$ over $\Q$ (corresponding to cyclic cubic fields)
is due to Cohn (see \cite{Cohn}), and over a general number field is due to 
the author and collaborators (see \cite{CoDiOl3}). The result of (1) over $\Q$
is certainly also in the literature (at least its main term), but over a 
general number field it seems to be new, as are all the other results, whether
over $\Q$ or over a general number field.

Note that the formula in (2) is given because of its elegance and for 
comparison with the quartic case, which we give below, but it should 
\emph{not} be used for practical computation of the constants 
$C(\Q(\sqrt{D}))$; for this, use instead Corollary \ref{corasym} below.
We emphasize that (for $D$ of reasonable size) all these constants can
easily be computed to hundreds of decimals, using the folklore method
explained in detail in Section 10.3.6 of \cite{Coh5}.

\subsection{Comparison with the Quartic Case}

Because of its striking similarity, we recall the results of \cite{Coh4}
in the \emph{quartic} case. Let $K_3$ be a cubic number field, and set 
$g(K_3)=3$ if $K_3$ is cyclic, $g(K_3)=1$ otherwise. We let $\mathcal{F}(K_3)$
be the set of isomorphism classes of quartic number fields $K$ whose cubic
resolvent is isomorphic to $K_3$. If $K\in\mathcal{F}(K_3)$, its discriminant 
$d(K)$ is of the form $d(K)=d(K_3)f(K)^2$ for some integer $f(K)$, and as in 
our case we let
$$M(K_3,X)=|\{K\in\mathcal{F}(K_3),\ f(K)\le X\}|\;.$$
The main result of \cite{Coh4} is then as follows:

\begin{theorem} Denote by $a_{K_3}(p)$ the number of copies of $\Q_p$ in 
$K_3\otimes\Q_p$ ($a_{K_3}(p)=0$, $1$ or $3$ according to whether the number 
of prime ideals above $p$ in $K_3$ is equal to $1$, $2$, or $3$). Then 
$M(K_3,X)=C(K_3)X+O(X^{\al})$ for some $\al<1$, with
$$C(K_3)=\dfrac{1}{g(K_3)}\dfrac{c_2(K_3)}{2^{4+r_2(K_3)}}\prod_{p\ne 2}\left(1+\dfrac{a_{K_3}(p)}{p}\right)\left(1-\dfrac{1}{p}\right)\;,$$
where
$$c_2(K_3)=\begin{cases}
11&\text{ if }2\Z_{K_3}=\p_1\\
14&\text{ if }2\Z_{K_3}=\p_1^3\\
15&\text{ if }2\Z_{K_3}=\p_1\p_2\\
16&\text{ if }2\Z_{K_3}=\p_1^2\p_2\text{ and }v_2(d(K_3))=3\\
18&\text{ if }2\Z_{K_3}=\p_1^2\p_2\text{ and }v_2(d(K_3))=2\\
23&\text{ if }2\Z_{K_3}=\p_1\p_2\p_3.\\
\end{cases}
$$
\end{theorem}
The similarities are striking.

\section{Galois Theory}

\begin{definition}We denote by $\rho=\z_3$ a primitive cube root of unity and
we set $L=K_2(\rho)$ and $k_z=k(\rho)$. We let $\tau$ be a generator of
$\Gal(L/K_2)$, and we let $\tau_2$
be a generator of $\Gal(K_2/k)$. We denote by
$G=\Gal(L/k)$. Finally, we let
$\sigma$ be one of the two generators of the cyclic group of order $3$
$\Gal(N/K_2)\isom\Gal(N_z/L)$, where $N_z=N(\rho)$.
\end{definition}

\noindent
{\bf Remark.}
We have the following relations:
$$\tau^2=\tau_2^2=1\;,\quad\tau\tau_2=\tau_2\tau\;,\quad\tau\sigma=\sigma\tau\;.$$

We will need to distinguish five cases, according to the triviality or not of
$\tau$ or $\tau_2$, and to their action on $\rho$. We will order them as
follows, and this numbering will be kept throughout the paper, so should be
referred to.

\begin{enumerate}\item $\tau=\tau_2=1$: here $K/k$ is a cyclic cubic extension;
in other words $K_2=k$, $\Gal(N_z/k)\isom C_3$, and $\rho\in k$.
\item $\tau_2=1$ and $\tau(\rho)=\rho^{-1}$: here $K/k$ is a cyclic cubic
extension, so that $K_2=k$, $\Gal(N_z/k)\isom C_6$; in other
words $\tau\sigma=\sigma\tau$, and $\rho\notin k$ so $L=k(\rho)$.
\item $\tau=1$ and $\tau_2(\rho)=\rho$ but $\tau_2\ne1$: here $K/k$ is
noncyclic, $\rho\in k$, and in particular $L=K_2$, and
$\Gal(N_z/k)\isom D_3$; in other words $\tau_2\sigma=\sigma^{-1}\tau_2$.
\item $\tau=1$ and $\tau_2(\rho)=\rho^{-1}$: here again $L=K_2$, so that $\rho\in K_2$, but $\rho\notin k$, so $K_2=k(\rho)$, and again
$\Gal(N_z/k)\isom D_3$; in other words $\tau_2\sigma=\sigma^{-1}\tau_2$.
\item $\tau\ne1$ and $\tau_2\ne1$: here $\rho\notin K_2$, so
$\tau(\rho)=\rho^{-1}$ but $\tau_2(\rho)=\rho$, so that the fixed field of
$L$ under $\tau_2$ is equal to $k_z=k(\rho)$, and
$\Gal(N_z/k)\isom D_3\times C_2$; in other words $\tau\sigma=\sigma\tau$ and
$\tau_2\sigma=\sigma^{-1}\tau_2$.
\end{enumerate}

\begin{definition}\begin{enumerate}\item In cases (1) to (5) above, we set
$T=\emptyset$, $\{\tau+1\}$, $\{\tau_2+1\}$, $\{\tau_2-1\}$, $\{\tau+1,\tau_2+1\}$,
respectively, where $T$ is considered as a subset of the group ring
$\Z[\Gal(L/k)]$ or of $\F_3[\Gal(L/k)]$.
\item We define $\i(\tau\pm1)=\tau\mp1$ and $\i(\tau_2\pm1)=\tau_2\mp1$.
\item For any group $M$ on which $T$ acts, we denote by $M[T]$ the subgroup
of elements of $M$ annihilated by all the elements of $T$.\end{enumerate}
\end{definition}
We will need the following trivial lemma (see \cite{CoDiOl3}, Lemma 2.4).

\begin{lemma}\label{lemit} Let $M$ be an $\F_3[G]$-module. For any $t\in T$ we
have $M[t]=\i(t)(M)$, and conversely $M[\i(t)]=t(M)$.\end{lemma}

\begin{proposition}\begin{enumerate}\item There exists a bijection between on
the one hand isomorphism classes of extensions $K/k$ having
quadratic resolvent field isomorphic to $K_2$, and on the other hand classes
of elements $\ov{\al}\in (L^*/{L^*}^3)[T]$ such that
$\ov{\al}\ne\ov{1}$ modulo the equivalence relation identifying $\ov{\al}$
with its inverse.
\item If $\al\in L^*$ is some representative of $\ov{\al}$,
the extension $K/k$ corresponding to $\al$ is the fixed field under
$\Gal(L/k)$ of the field $N_z=L(\root3\of{\al})$.\end{enumerate}
\end{proposition}

\Proof Since $\rho\in L$, by Kummer theory, cyclic cubic extensions of $L$ are of the form $N_z=L(\root3\of\al)$, where
$\ov{\al}\ne\ov{1}$ is unique in $(L^*/{L^*}^3)$
modulo the equivalence relation identifying $\ov{\al}$ with its inverse.
If $\th^3=\al$, then we may
assume that $\sigma(\th)=\rho\th$. When $\tau$ is non trivial (cases
(2) and (5)) we have $\tau(\rho)=\rho^{-1}$. Thus,
$$\sigma(\th\tau(\th))=\rho\th\tau(\sigma(\th))=\rho\th\tau(\rho\th)=\th\tau(\th)\;,$$
so by Galois theory $\th\tau(\th)\in L$, so
$\al\tau(\al)$ is a cube, in other words
$\al\in (L^*/{L^*}^3)[\tau+1]$.

Similarly when  $\tau_2$ is nontrivial we
have either $\tau_2(\rho)=\rho$
(cases (3) and (5)) or  $\tau_2(\rho)=\rho^{-1}$ (case (4)). A similar
computation gives respectively $\al\in (L^*/{L^*}^3)[\tau_2+1]$ or
$\al\in (L^*/{L^*}^3)[\tau_2-1]$.

Conversely, assume that these conditions are satisfied. The group
conditions on $\tau$ and $\tau_2$ are automatically satisfied, and the group conditions on $\sigma$ are exactly those
corresponding to the set $T$. It follows that $N_z/k$ is Galois with suitable
Galois group. The uniqueness statement comes from the corresponding statement
of Kummer theory, since $\al$ and $\al^{-1}$ give the same extension.\fp

\smallskip

\begin{definition} We denote by $V_3(L)$ the group of
$(3-)$\emph{virtual units} of $L$, in other words the group of
$u\in L^*$ such that $u\Z_L=\q^3$ for some ideal $\q$ of
$L$. We define the $(3-)$\emph{Selmer group} $S_3(L)$ of
$L$ by $S_3(L)=V_3(L)/{L^*}^3$.\end{definition}

It is immediate that the Selmer group is finite.

\begin{proposition}\label{propbij}
\begin{enumerate}\item There exists a bijection between isomorphism classes of
cubic extensions $K/k$ with given quadratic resolvent field $K_2$
and equivalence classes of triples $(\a_0,\a_1,\ov{u})$ modulo the equivalence
relation $(\a_0,\a_1,\ov{u})\sim(\a_1,\a_0,1/\ov{u})$, where $\a_0$, $\a_1$,
and $\ov{u}$ are as follows:
\begin{enumerate}\item The $\a_i$ are coprime integral squarefree ideals of $L$
such that $\overline{\a_0\a_1^2}\in Cl(L)^3$ and $\a_0\a_1^2\in (I/I^3)[T]$, 
where $I$ is the group of fractionals ideals of $L$.
\item $\ov{u}\in S_3(L)[T]$, and $\ov{u}\ne1$ when
$\a_0=\a_1=\Z_L$.
\end{enumerate}
\item If $(\a_0,\a_1)$ is a pair of ideals satisfying (a) there exist an ideal
$\q_0$ and an element $\al_0$ of $L$ such that $\a_0\a_1^2\q_0^3=\al_0\Z_L$ 
with $\al_0\in (L^*/{L^*}^3)[T]$. The cubic extensions $K/k$ corresponding to
such a pair $(\a_0,\a_1)$ are given as follows: for any $\ov{u}\in S_3(L)[T]$
the extension is the cubic subextension of $N_z=L(\root3\of{\al_0u})$
(for any lift $u$ of $\ov{u}$).
\end{enumerate}
\end{proposition}

\Proof Let $N_z=L(\root3\of\al)$ as above. We can write uniquely
$\al\Z_L=\a_0\a_1^2\q^3$ where the $\a_i$ are coprime squarefree
ideals of $L$. Since $\al\in
(L^*/{L^*}^3)[T]$ and the class of $\a_0\a_1^2$ is equal to that of
$\q^{-3}$, we obtain  \emph{(a)}.
Now let $\a_0,\a_1$ be given satisfying \emph{(a)}. There exists an ideal
$\q$ and an element
$\al\in L$ such that $(\a_0\a_1^2)\q^3=\al\Z_L$. Applying
any $t\in T$, we deduce that
$\q_1^3=t(\al)\Z_L$ for some ideal $\q_1$, so that $t(\al)$ is a virtual 
unit.
From $t\circ \i(t)=0$ and
Lemma \ref{lemit} we deduce that $t(\al)\in t(S_3(L))$, in other
words that $t(\al)=\ga^3 t(u)$,
for some virtual unit $u$ and some element $\ga$. Thus, if we set
$\al_0=\al/u$, we have $\al_0\in (L^*/{L^*}^3)[t]$, and
$\a_0\a_1^2\q_0^3=\al_0\Z_L$, for some ideal $\q_0$.

The rest of the proof is immediate: $\a_0\a_1^2\q_0^3=\al_0\Z_L$
and $\a_0\a_1^2\q^3=\al\Z_L$,
with both $\al_0,\al\in (L^*/{L^*}^3)[T]$
if and only if $\al/\al_0=(\q/\q_0)^3\in V_3(L)[T]$, so $\al=\al_0u$ for
some lift $u$ of $\ov{u}\in S_3(L)[T]$. Finally $\al$ and $\be$ give
equivalent extensions if and only if either $\be=\al\ga^3$, which does not
change the $\a_i$ and the class $\ov{u}$,
or if $\be=\al^{-1}\ga^3$. In this case
$$\be\Z_L=\a_0^{-1}\a_1^{-2}\q^{-3}\ga^3=\a_1\a_0^2(\ga\a_0^{-1}\a_1^{-1}\q^{-1})^3\;,$$
which interchanges $\a_0$ and $\a_1$, and changes $\ov{u}$ into $1/\ov{u}$, finishing the proof. Note
that the only fixed point of this involution on triples is obtained for
 $\a_0=\a_1=\Z_L$, and $\ov{u}=1$.\fp

\begin{lemma}\label{lema0a1}
\begin{enumerate}
\item  The condition $\a_0\a_1^2\in(I/I^3)[T]$ is equivalent to 
$\a_1=\tau(\a_0)$, $\a_1=\tau_2(\a_0)$, $\a_0=\tau_2(\a_0)$ and 
$\a_1=\tau_2(\a_1)$, and $\a_1=\tau(\a_0)=\tau_2(\a_0)$ in cases (2), (3), 
(4), and (5), respectively.
\item The ideal $\a_0\a_1$ of $L$ comes from an ideal $\a_{\al}$ of $K_2$
(in other words $\a_0\a_1=\a_{\al}\Z_L$), and in cases (1), (2), and (3)
it comes from an ideal of $k$, while in cases (4) and (5), $\a_{\al}$ is an
ideal of $K_2$ invariant by $\tau_2$.\end{enumerate}
\end{lemma}

\Proof 
Just apply uniqueness of decomposition to $\tau(\a_0\a_1^2)$ and $\tau_2(\a_0\a_1^2)$.\fp

In case (5), which is the only case where $G=\Gal(L/k)\isom V_4$, we define 
$K_2'$ to be the quadratic subextension of $L/k$ different from $K_2$ and 
$k_z$. 
\begin{definition}\label{defD}
We define $\mathcal{D}$ (resp., $\mathcal{D}_3$) to be the set of all 
prime ideals $p$ in $k$ with $p\nmid 3\Z_k$ (resp., with $p\mid 3\Z_k$), such 
that:
\begin{itemize}
\item no other conditions  in cases (1) and (4);
\item $p$ is split in $L/k$ in case (2) and (3);
\item the ideals above $p$ are split in $L/K_2$ and $L/k_z$ in case (5).
\end{itemize}
\end{definition}

\begin{proposition}\label{split}
\begin{enumerate}\item Let $\p$ be a prime ideal of $K_2$ dividing $\a_{\al}$
and let $p$ be the prime ideal of $k$ below $\p$. Then
$p\in\mathcal{D}\cup \mathcal{D}_3$.
\item In cases (2) and (3), set $K'_2=L$. Then in cases (2), (3), and (5) we 
have $p\in\mathcal{D}\cup \mathcal{D}_3$ if and only $p$ is split in 
$K'_2/k$.\end{enumerate}
\end{proposition}

\Proof (1) is immediate, for (2) use decomposition groups.\fp

\section{Conductors}

The discriminant (equivalently, the conductor) of a cyclic Kummer
extension is given by an important theorem of Hecke (see \cite{Coh2}, Section
10.2.9). We will mainly need it in the cubic case, but we also need it in
the quadratic case, where it takes an especially nice form:

\begin{theorem}\label{thk2} Let $k$ be a number field, let $K_2=k(\sqrt{D})$
be a quadratic extension with $D\in k^*\setminus {k^*}^2$, and write uniquely
$D\Z_k=\a\q^2$, where $\a$ is an integral squarefree ideal. Then
$$\gd(K_2/k)=\f(K_2/k)=4\a/\c^2\;,$$
where $\c$ is the largest ideal (for divisibility) dividing $2\Z_k$ and
coprime to $\a$ such that the congruence $x^2/D\equiv 1\pmod{^*\c^2}$ has a
solution.\end{theorem}

\begin{corollary}\label{ramK2} Let $K$ be a number field such that
$\rho\notin k$, where $\rho=\z_3$ is a primitive cube root of unity,
and set $K_z=K(\rho)$. Then
$$\gd(K_z/K)=\prod_{\substack{\p\mid 3\Z_K\\e(\p/3)\text{ odd}}}\p\;.$$
In particular, the ramified primes in $K_z/K$ are those above $3$ such that
$e(\p/3)$ is odd.
\end{corollary}

\Proof We have $K_z=K(\sqrt{-3})$, so $D=-3$. We have
$D\Z_K=3\Z_K=\a\q^2$ with
$\a=\prod_{\substack{\p\mid 3\Z_K\\e(\p/3)\text{ odd}}}\p\;.$
On the other hand $\a$ is coprime to $2$ and the congruence
$x^2\equiv-3\pmod{4}$ has the solution $x=1$, so $\c=2\Z_K$ and the corollary
follows.\fp
If $\p$ is a prime ideal of $K_2$, we will denote by $\p_z$ any prime ideal
of $L$ above $\p$. By the above corollary, we have $e(\p_z/\p)=2$
if and only if $L\ne K_2$ and $e(\p/3)$ is odd, otherwise
$e(\p_z/\p)=1$.

\vskip.5cm

In the case of cyclic cubic extensions, the result is more complicated,
especially when $L\ne K_2$. We first need some definitions.

\begin{definition}\label{definitionot} In the sequel, when $p$ is a prime ideal of 
$k$ we will denote by $\p$ a prime ideal of $K_2$ above $p$, and by
$\p_z$ a prime ideal of $L$ above $\p$. In addition, to simplify notation:
\begin{itemize}
\item We set $p^{1/2}=\p$ if $p$ is ramified in $K_2/k$ (i.e., 
$p\Z_{K_2}=\p^2$), and similarly $\p^{1/2}=\p_z$ if $\p$ is ramified in 
$L/K_2$ (i.e., $\p\Z_L=\p_z^2$).
\item We say that $p\subset k$ divides some ideal $\b$ of $K_2$ (resp., of $L$)
when $(p\Z_{K_2})^{1/e(\p/p)}$ (resp., $(p\Z_L)^{1/e(\p_z/p)}$) does. 
\end{itemize}
\end{definition}

Note that $e(\p_z/p)\le2$ (indeed, if for instance $e(\p/p)=2$ then $e(\p/3)$
is even so $\p_z/\p$ is unramified by Corollary \ref{ramK2}), so we will never
need to define ``$p^{1/4}$''.

\begin{definition}\label{defaal1} Let $\ov{\al}\in(L^*/{L^*}^3)[T]$ be as 
above, let $p$ be an ideal of $k$ above $3$, let $\p$ and $\p_z$ be as in
Definition \ref{definitionot}, and consider the congruence 
$x^3/\al\equiv1\pmod{^*\p_z^n}$ in $L$. If this congruence is soluble for 
$n=3e(\p_z/3)/2$ we set $A_{\al}(p)=3e(\p_z/3)/2+1$, otherwise, if
$n<3e(\p_z/3)/2$ is the largest exponent for which it has a solution, we set
$A_{\al}(p)=n$. In both cases we define
$$a_{\al}(p)=\dfrac{A_\al(p)-1}{e(\p_z/p)}\;.$$
\end{definition}
It is clear that $A_{\al}(p)$ and $a_{\al}(p)$ do not depend on the
ideal $\p_z$ above $p$, whence the notation. We have the following properties:

\begin{proposition}\label{propaal} We have
$0\le a_{\al}(p)< 3e(p/3)/2-1/e(\p/p)$ and $a_{\al}(p) e(\p/p)\in\Z$, or
$a_{\al}(p)=3e(p/3)/2$, which happens if and only if
$A_{\al}(p)=3e(\p_z/3)/2+1$, in which case it is only a half integer
when $e(\p_z/\p)=2$.
\end{proposition}

\begin{definition}
To simplify notations, we set 
$$\mathcal{P}_3=\{p\mid 3\mathbb{Z}_k\textrm{ such that } e(\p/3) \textrm{ odd}\}\;.$$
\end{definition}

\begin{theorem}\label{condN} Let $N$ correspond to $\al$ as above, write
uniquely $\al\Z_L=\a_0\a_1^2\q^3$ with $\a_0$ and $\a_1$ integral coprime
squarefree ideals, and let $\a_{\al}$ be the ideal of $K_2$ such that
$\a_0\a_1=\a_{\al}\Z_L$ (see Lemma \ref{lema0a1}). Then
$$\f(N/K_2)=\dfrac{3\a_{\al}\prod_{p\mid 3\Z_{k}}(p\Z_{K_2})^{e(p/3)/2}
\prod_{p\in\mathcal{P}_3} (p\Z_{K_2})^{1/2}}
{\prod_{\substack{p\mid3\Z_{k}\\p\nmid\a_{\al}}}(p\Z_{K_2})^{\lceil a_{\al}(p) e(\p/p)\rceil/e(\p/p)}}\;.$$
\end{theorem}
{\bf Remark.}  Proposition \ref{propaal} and Theorem \ref{condN} come from similar results in \cite{CoDiOl3}
where we have just replaced $a_{\al}(\p)$ by $a_{\al}(p)=a_{\al}(\p)/e(\p/p)$.
In particular, the fact that $a_{\al}(p)e(\p/p)$ is an integer when 
$a_{\al}(p)<3e(p/3)/2$ is a rather subtle result, which follows from the 
use of higher ramification groups. 

\begin{definition}\label{defh} Let $p$, $\p$ and $\p_z$ be as in Definition \ref{defaal1},
and let $a$ be such that $0\le a<3e(p/3)/2-1/e(\p/p)$ and $a e(\p/p) \in\Z$, 
or $a=3e(p/3)/2$. For $\eps=0$ or $1$ we define $h(\eps,a,p)$ as follows:
\begin{itemize}
\item We set $h(0,a,p)=0$ if  $a=3e(p/3)/2$ or  $e(\p_z/\p)=2$; in the
  other cases we set $h(0,a,p)=1/e(\p/p)$.
\item we set $h(1,a,p)=2/e(\p_z/p)$.
\end{itemize}\end{definition}

\begin{lemma}\label{lembiai} Let $b=a+h(\eps,a,p)$.
\begin{enumerate}\item Assume that $b\le 3e(p/3)/2$. Then
$h(\eps,b,p)=h(\eps,a,p)$, so that $a=b-h(\eps,b,p)$.
\item We have $b=0$ if and only if $a=0$, $\eps=0$, and $e(\p_z/\p)=2$.
In particular, if $e(\p_z/\p)=1$ we have $b>0$.
\end{enumerate}\end{lemma}

\Proof Follows immediately from Definition \ref{defh}.\fp

\begin{lemma}\label{lembdef} Let $p$ be a prime ideal of $k$ and denote
by $D_k$ the congruence $x^3/\al\equiv1\pmod{^*p^k}$ in $L$. If $a$ is as in
the above definition, then $a_{\al}(p)=a$ if and only if $D_k$ is soluble for
$k=a+h(0,a,p)$ and not soluble for $k=a+h(1,a,p)$, where this last
condition is ignored if $a+h(1,a,p)>3e(p/3)/2$.\end{lemma}

\Proof Just apply definitions \ref{defaal1} and \ref{defh} and
Proposition \ref{propaal}.\fp

\section{The Dirichlet Series}\label{sec4}

To avoid having both the norm from $K_2/\Q$ and from $k/\Q$, and to
emphasize the fact that we are mainly interested in the latter, we set
explicitly the following definition:

\begin{definition} If $\a$ is an ideal of $k$, we set $\N(\a)=\N_{k/\Q}(\a)$,
while if $\a$ is an ideal of $K_2$, we set
$$\N(\a)=\N_{K_2/\Q}(\a)^{1/[K_2:k]}\;.$$
\end{definition}

This practical abuse of notation cannot create any problems since
if $\a$ is an ideal of $k$ we have $\N(\a)=\N(\a\Z_{K_2})$.
For instance, since $\f(N/K_2)=\f(K/k)\Z_{K_2}$, we have 
$\N(\f(K/k))=\N(\f(N/K_2))$. We emphasize that unless explicitly
written otherwise, from now on we will only use the above notation.

\begin{definition} The fundamental Dirichlet series is defined by
$$\Phi(s)=\dfrac{1}{2}+\sum_{K\in\mathcal{F}(K_2)}\dfrac{1}{\N(\f(K/k))^s}\;.$$
\end{definition}

\begin{definition} For $\al_0\in L^*$ and $\b$ an ideal of $L$ we introduce 
the function
$$f_{\al_0}(\b)=\bigl|\{\overline{u}\in S_3(L)[T],\ x^3/(\al_0u)\equiv 1 \pmod {^*\b}\text{\quad soluble in } L\}\bigr|\;,$$
with the convention that $f_{\al_0}(\b)=0$ if $\b\nmid 3\sqrt{-3}.$
\end{definition}

\begin{definition}\label{defB}
\begin{enumerate}\item We let $\mathcal{B}$ be the set of formal products of 
the form \\${\prod_{p_i\mid 3\Z_k}(p_i\Z_{K_2})^{b_i}}$, where the $b_i$ are
such that $0\le b_i\le 3e(p_i/3)/2$ and 
$e(\p_i/p_i)b_i\in\Z\cup \{3e(\p_i/3)/2\}$.
\item We will consider any $\b\in\mathcal{B}$ as an ideal of $K_2$, where by
abuse of language we accept to have half powers of prime ideals of $K_2$, and
we set $\b_z=\b\Z_L$.
\item If $\b=\prod_{\p_i\mid 3\Z_{K_2}}\p_i^{{b_i}'}\in\mathcal{B}$, ${b_i}'=e(\p_i/p_i)b_i$, we set
$\NC(\b)=\prod_{\p_i\mid\b}\N(\p_i)^{\lceil {b_i}'\rceil}\;.$
\item For $\b\in\mathcal{B}$ we define
  $\gothr^e(\b)=\prod_{\substack{\p\mid 3\Z_{K_2},\ \p\nmid
    \b\\ e(\p/3) \text{ even}}} \p$.
\item We set $\gd_3=\prod_{p\in\mathcal{D}_3}p$.
\end{enumerate}
\end{definition}

\begin{definition}\label{lempb}
\begin{enumerate}
\item Set $e=e(p/3)$, let
$\p$ an ideal of $K_2$ above $p$, let $\p_z$ be an ideal of $L$ above $\p$,
and define $s'=s/e(\p/p)$. We define $Q((p\mathbb{Z}_{K_2 })^{b_i},s)$
as follows:
\begin{itemize}
\item if $e(\p_z/\p)=1$, (so that $e(\p/3)$ is even) we have 
$$\kern-10pt Q((p\Z_{K_2})^b,s)=\begin{cases} 
0 & \text{\kern-6pt if } b=0\;,\\
1/\N(p)^{s'} & \text{\kern-6pt if }b=1/e(\p/p)\;,\\
1/\N(p)^{s'}-1/\N(p)^{2s'} & \text{\kern-6pt if } 2/e(\p/p)\le b \le 3e/2-1/e(\p/p)\;,\\
1-1/\N(p)^{2s'} & \text{\kern-6pt if } b=3e/2\;.\end{cases}$$
\item if $e(\p_z/\p)=2$ (so that $e(\p/p)=1$) we have
$$Q((p\Z_{K_2})^b,s)=\begin{cases}
1 & \text{ if } b=0 \textrm{ or }b=3e/2\;,\\
1-1/\N(p)^{s'} & \text{ if } 1\le b\le 3e/2-3/2\;,\\
-1/\N(p)^{s'} & \text{ if } b=3e/2-1/2\;. \end{cases}$$
 \end{itemize}
\item We set $P_\b(s)=\prod_{p\mid \b}
Q((p\Z_{K_2})^{v_{p}(\b)},s)$.
\end{enumerate}
\end{definition}

\begin{proposition}
We have
$$ \Phi(s)=\dfrac{1}{2\cdot 3^{(3/2)[k:\Q]s}\prod_{\mathcal{P}_3}\N(p)^{s/2}}\sum_{\substack{\b\in\mathcal{B}\\\gothr^e(\b)\mid\gd_3}} \NC (\b)^s
  P_\b(s)\sum_{\substack{(\a_0,\a_1)\in J
 \\ (\a_{\al},3\Z_{K_2})=\gothr^e(\b)}}
\dfrac{f_{\al_0}(\b)}{\N(\a_{\al})^s}\;.$$
\end{proposition}
\Proof
This formula is obtained after some computations, applying in
particular Proposition \ref{propbij}, Theorem \ref{condN} and an
inclusion-exclusion argument. A complete proof of the analogous result
in the (simpler) case of cyclic extensions can be found in \cite{CoDiOl3}.
\fp

\section{Computation of $f_{\al_0}(\b)$}\label{sec5}

Recall that $\b_z\mid3\sqrt{-3}$ and that the $\a_i$ are coprime squarefree
ideals such that $\a_0\a_1^2\in (I/I^3)[T]$ and 
$\overline{\a_0\a_1^2}\in Cl(L)^3$. We have also set
$\a_0\a_1^2\q_0^3=\al_0\Z_L$ with $\al_0\in (L^*/{L^*}^3)[T]$. Recall that
$$f_{\al_0}(\b)=\bigl|\{\overline{u}\in S_3(L)[T],\ x^3\equiv\al_0u\pmod {^*\b_z}\text{\quad soluble in } L\}\bigr|\;,$$
where we have replaced the congruence $x^3/(\al_0u)\equiv1\pmod{^*\b_z}$ by
the above since we may assume $\al_0$ coprime to $\b_z$ (changing
$\q_0$ and $\al_0$ if necessary).

\begin{definition}\label{sbk} Set
$$S_{\b}(L)[T]=\{\ov{u}\in S_3(L)[T],\ x^3\equiv u\pmod{^*\b_z}\text{ soluble}\}\;,$$
where $u$ is any lift of $\ov{u}$ coprime to $\b_z$, and the congruence is
in $L$.
\end{definition}

\begin{lemma}\label{lemfal1et2}
Let $\a_0,\a_1$ as in condition (1) of Proposition \ref{propbij}. Then
$$f_{\al_0}(\b)=\begin{cases}|S_\b(L)[T]|& \text{if
  $\overline{\a_0\a_1^2}\in Cl_\b(L)^3$}\\ 0 & \text{otherwise.}\end{cases}$$
\end{lemma}
\Proof
First, assume that there exists an $u_0\in S_3(L)[T]$ such that 
$x_0^3\equiv\al_0u_0\pmod{^*\b_z}$ for some $x_0\in L$. The congruence 
$x^3\equiv\al_0u\pmod{^*\b_z}$ is thus equivalent to
$(x/x_0)^3\equiv(u/u_0)\pmod{^*\b_z}$, in other words to 
$u/u_0\in S_{\b}(L)[T]$, so the set of possible $\ov{u}$ is equal to 
$\ov{u_0}S_{\b}(L)[T]$, whose cardinality is $|S_{\b}(L)[T]|$.
So if $f_{\al_0}\neq 0$ then it is equal to $|S_\b(L)[T]|$.

Now let us prove that  $f_{\al_0}\neq 0$ if and only if  $\overline{\a_0\a_1^2}\in Cl_\b(L)^3$. The condition $\overline{\a_0\a_1^2}\in Cl_{\b}(L)^3$ is equivalent to
the existence of $\q_1$ and $\be_1\equiv1\pmod{^*\b_z}$ such that
$\a_0\a_1^2\q_1^3=\be_1\Z_L$. Assume first that $u$ exists, so that
$x_0^3=\al_0u\be$ for some $\be\equiv1\pmod{^*\b_z}$ and $u\Z_L=\q^3$. It
follows that $\a_0\a_1^2\q_0^3\q^3=\al_0u\Z_L=(x_0^3/\be)\Z_L$,
so we can take $\q_1=\q_0\q/x_0$ and $\be_1=1/\be\equiv1\pmod{^*\b_z}$.
Conversely, assume that $\a_0\a_1^2\q_1^3=\be_1\Z_L$ with
$\be_1\equiv1\pmod{^*\b_z}$.  Since $\a_0\a_1^2\in(I/I^3)[T]$, we have
$t(\be_1)=\ga^3$ for some $\ga\in L^*$.
It follows that $\al_0\Z_L=\a_0\a_1^2\q_0^3=\be_1(\q_0/\q_1)^3$. Thus,
$u=\al_0/\be_1$ is a virtual unit, and $t(u)$ is a cube of $L$
since this is true for $\al_0$ and for $\be_1$. Thus 
$\ov{u}\in S_3(L)[T]$ and $1^3\equiv\be_1\equiv\al_0/u\pmod{^*\b_z}$,
so $f_{\al_0}(\b)\ne0$, proving the lemma.\fp

Note that when we assume $\overline{\a_0\a_1^2}\in Cl_\b(L)^3$ we have
automatically $\overline{\a_0\a_1^2}\in Cl(L)^3$, so we only need to
assume that $\a_0\a_1^2\in (I/I^3)[T]$.

To compute  $|S_{\b}(L)[T]|$ we will use
the folling lemmas, which are similar to the ones proposed in
(\cite{CoDiOl3}, \S 2), so we will omit the proofs.

\begin{lemma}\label{lemfal3} Set $Z_{\b}=(\Z_L/\b_z)^*$, $Cl=Cl(L)$,
  $Cl_\b=Cl_\b(L)$ and $U=U(L)$. Then
$$|S_{\b}(L)[T]|=\dfrac{|(U/U^3)[T]| |(Cl_{\b}/Cl_{\b}^3)[T]|}{|(Z_{\b}/Z_{\b}^3)[T]|}\;.$$
In particular
$$|S_3(L)[T]|=|(U/U^3)[T]||(Cl/Cl^3)[T]|\;.$$
\end{lemma}

The quantity $|(Cl_{\b}/Cl_{\b}^3)[T]|$ will in fact disappear
in subsequent computations, and in any case cannot be computed more
explicitly. 

\begin{lemma}\label{lemu5} For any number field $K$, denote by $\rk_3(K)$
the $3$-rank of the group of units of $K$, in other words
$\rk_3(K)=\dim_{\F_3}(U(K)/U(K)^3)$, so that $|U(K)/U(K)^3|=3^{\rk_3(K)}$.
\begin{enumerate}\item With evident notation we have
$$\rk_3(K)=\begin{cases} 
r_1(K)+r_2(K)-1 & \text{ if $\rho\notin K$,}\\
r_1(K)+r_2(K) & \text{ if $\rho\in K$.}\end{cases}$$
\item We have $|(U/U^3)[T]|=3^{r(U)}$, where
$$r(U)=\begin{cases}
\rk_3(k) & \text{ in cases (1) and (4),}\\
\rk_3(L)-\rk_3(k) & \text{ in cases (2) and (3),}\\
\rk_3(L)+\rk_3(k)-\rk_3(K_2)-\rk_3(k_z) & \text{ in case (5).}\end{cases}$$
\end{enumerate}
\end{lemma}

\begin{lemma}\label{lemfal4} Assume that $\b$ is an ideal of $\mathcal{B}$,
stable by $\tau_2$ and such that ${\b_z\mid3\sqrt{-3}}$, and define
$$\c_z=\prod_{\substack{\p_z\subset L\\ \p_z\mid\b_z}}\p_z^{\lceil v_{\p_z}(\b_z)/3\rceil}\;.$$
Then
$$|(Z_{\b}/Z_{\b}^3)[T]|=\left|\dfrac{\c_z}{\b_z}[T]\right|\;.$$
\end{lemma}

\begin{lemma}\label{lemcb}
$$|(Z_{\b}/Z_{\b}^3)[T]|=\begin{cases} |\c_z/\b_z| &\text{ in case (1)}\\
    \dfrac{|\c_z/\b_z|}{|(\c_z\cap k)/(\b_z\cap k)|} &\text{\kern-38pt in cases (2) and
    (3)}\\ |(\c_z\cap k)/(\b_z\cap k)| & \text{in case (4)}\\
    \dfrac{|\c_z/\b_z| |(\c_z\cap k)/(\b_z\cap k)|}{|(\c_z\cap K_2)/(\b_z\cap K_2)|
    |(\b_z\cap k_z)/(\c_z\cap k_z)|}& \text{in case (5).}\end{cases}$$
\end{lemma}

\section{Final Form of the Dirichlet Series}

We can now put together all the work that we have done. Recall that we
have computed $|U/U^3[T]|$ in Lemma \ref{lemu5} and
$|(Z_\b/Z_\b^3)[T]|$ in Lemma \ref{lemcb}. Moreover,
$\mathcal{B}$, $\NC$ and $\gd_3$ are defined in Definition
\ref{defB}, and $P_\b(s)$
is given by Definition \ref{lempb}. Finally, recall that
we have
$$\Phi(s)=\dfrac{1}{2}+\sum_{K\in\mathcal{F}(K_2)}\dfrac{1}{\N(\f(K/k))^s}\;.$$

\begin{theorem}\label{mainth1} For any ideal $\b$, set 
$G_{\b}=(Cl_{\b}/Cl_{\b}^3)[T]$. We have
\begin{align*}\Phi(s)&=\dfrac{|(U/U^3)[T]|}{2\cdot 3^{(3/2)[k:\Q]s}
\prod_{\mathcal{P}_3}\N(p)^{s/2}}\sum_{\substack{\b\in\mathcal{B}\\\gothr^e(\b)\mid\gd_3}}
\left(\dfrac{\NC(\b)}{\N(\gothr^e(\b))}\right)^s\dfrac{P_{\b}(s)}{|(Z_{\b}/Z_{\b}^3)[T]|}\sum_{\chi\in\widehat{G_{\b}}}F(\b,\chi,s) \;,\end{align*}
where
$$F(\b,\chi,s)=\prod_{\substack{p\mid\gothr^e(\b)\\p\in\mathcal{D}_3'(\chi)}}2\prod_{\substack{p\mid\gothr^e(\b)\\p\in\mathcal{D}_3\setminus{\mathcal{D}_3}'(\chi)}}(-1)\prod_{p\in\mathcal{D}'(\chi)} \left(1+\dfrac{2}{\N(p)^{s}}\right)\prod_{p\in\mathcal{D}\setminus\mathcal{D}'(\chi)} \left(1-\dfrac{1}{\N(p)^{s}}\right)\;,$$
and $\mathcal{D}'(\chi)$ (resp. $\mathcal{D}_3'(\chi)$) is the set of
$p\in\mathcal{D}$ (resp. $\mathcal{D}_3$) such that $\chi(p\Z_L)=1$ in
cases (1) and (4), $\chi(\c)=\chi(\tau'(\c))$ in the other cases,
where we write $p\Z_L=\c\tau'(\c)$, $\tau'\in\{\tau,\tau_2\}$, and
$\c$ is not necessarily a prime ideal.

\end{theorem}

\Proof Let $\a_0$ and $\a_1$ be
as in condition (a) of Proposition \ref{propbij}. We
have $\a_0\a_1^2 \in Cl_{\b}(L)^3$ if and only if $\chi(\a_0\a_1^2)=1$ for all
characters $\chi\in\widehat{G_{\b}}$. The number of such characters being 
equal to $|G_{\b}|$, by orthogonality of characters we have

$$\Phi(s)=\dfrac{|(U/U^3)[T]|}{2\cdot 3^{(3/2)[k:\Q]s}
\prod_{\substack{p\mid 3\Z_{k},\\ e(\p/3)\text{	odd}}}\N(p)^{s/2}}\sum_{\substack{\b\in\mathcal{B}\\\gothr^e(\b)\mid\gd_3}}\dfrac{ \NC(\b)^s P_{\b}(s)}{|(Z_{\b}/Z_{\b}^3)[T]|}\sum_{\chi\in\widehat{G_{\b}}}H(\b,\chi,s)\;,$$
with
$H(\b,\chi,s)=\sum_{\substack{(\a_0,\a_1)\in J'\\ (\a_{\al},3\Z_{K_2})=\gothr^e(\b)}}\dfrac{\chi(\a_0\a_1^2)}{\N(\a_{\al})^s}\;,$
where $J'$ is the set of pairs of coprime squarefree ideals of 
$L$, satisfying condition (1) of Lemma \ref{lema0a1}, with no class group 
condition. Thus
$$H(\b,\chi,s)=\dfrac{\chi(\gothr^e(\b))}{\N(\gothr^e(\b))^s}\sum_{\substack{(\a,3\Z_L)=1\\\a \text{ squarefree }\\ \tau(\a)=\tau_2(\a)=\a}}\dfrac{\chi(\a)}{\N(\a)^s}\sum_{\a_1\mid\a\gothr^e(\b),\ \a_1\in J''}\chi(\a_1)\;,$$
where  $J''$ is the set of squarefree ideals $\a_1$ such that $\a_1$ is stable
by $\tau_2$ in case (4), $\a_1\tau'(\a_1)=\a\gothr^e(\b)$ for each nontrivial
$\tau'\in \{\tau,\tau_2\}$ in the other cases.

Let us define $G(\chi,p)$ by:
$$G(\chi,p)=\begin{cases}1+\chi(p\Z_L)& \text{ in cases (1) and (4), and
 otherwise :} \\ \chi(\c)+\chi(\tau'(\c)) & \text{ when
 }p\Z_L=\c\tau'(\c),\text{($\tau'$ and $\c$ as above)}.\end{cases}$$

Since $\a$ is coprime to $3$, by multiplicativity we have
$H(\b,\chi,s)=S_1S_2$ with
\begin{align*}S_1&=\dfrac{\chi(\gothr^e(\b))}{\N(\gothr^e(\b))^s}\prod_{p\mid\gothr^e(\b)}G(\chi,p)\text{\quad and\quad}\\
S_2&=\sum_{\substack{(\a,3\Z_L)=1\\\a\text{ squarefree }\\\tau(\a)=\tau_2(\a)=\a}}\dfrac{\chi(\a)}{\N(\a)^s}\prod_{p\mid\a}G(\chi,p)=
\prod_{p\in\mathcal{D}}\left(1+\dfrac{\chi(p\Z_L)G(\chi,p)}{\N(p)^{s}}\right)\;.\end{align*}
Now, looking at
the possible values for $G(\chi,p)$, we conclude.\fp

\begin{corollary}\label{resgen1} In cases (2) and (3), set $K'_2=L$, and in
all cases
denote by $\gd(K'_2/k)$ the relative discriminant of $K'_2/k$. Let us define
$$c_1=\dfrac{|(U/U^3)[T]|}{2\cdot3^{(3/2)[k:\Q]}\prod_{\substack{p\mid 3\Z_{k}\\e(\p/3) \text{ odd}}} \N(p)^{1/2}}\;,$$
$$c_2=\sum_{\substack{\b\in\mathcal{B}\\\gothr^e(\b)\mid\gd_3}}\dfrac{\NC (\b)}{\N(\gothr^e(\b))}\dfrac{P_{\b}(1)}{|(Z_{\b}/Z_{\b}^3) [T]|}2^{\om(\gothr^e(\b))}\;,$$
$$c_3=\prod_{p\subset k}\left(1-\dfrac{3}{\N(p)^{2}}+\dfrac{2}{\N(p)^{3}}\right)\prod_{p\mid 3\Z_k}\left(1+\dfrac{2}{\N(p)}\right)^{-1}\;,$$
$$c_4=\dfrac{1}{\z_k(2)}\prod_{p\in\mathcal{D}}
  \left(1-\dfrac{2}{\N(p)(\N(p)+1)}\right)\prod_{p\mid\gd(K_2'/k)}\left(1-\dfrac{1}{\N(p)+1}\right)\;,$$
where $\om(\gothr^e(\b))=\sum_{p\mid\gothr^e(\b)}1$.
\begin{itemize}
\item In cases (1) and (4), around $s=1$ we have
$$\Phi(s)=\dfrac{C(K_2/k)}{(s-1)^2}+\dfrac{C(K_2/k)D(K_2/k)}{s-1}+O(1)\;,$$
with constants 
\begin{align*}
C(K_2/k)&=c_1c_2c_3(\Res_{s=1}\z_k(s))^2\text{\quad and}\\
D(K_2/k)&=2\ga_k+\lim_{s\to1}\dfrac{G'(s)}{G(s)}\,\text{ where}\\
G(s)&=\dfrac{\Phi(s)}{\z_k(s)^2}\text{\quad and\quad}\ga_k=\lim_{s\to1}\left(\dfrac{\z_k(s)}{\Res_{s=1}\z_k(s)}-\dfrac{1}{s-1}\right)\;.\end{align*}

In addition, using the notation given at the beginning of this paper,
as $X\to\infty$ we have 
$$M(K_2/k,X)=C(K_2/k)X(\log(X)+D(K_2/k)-1)+O(X^{\alpha})\text{\quad for some $\al<1$\;.}$$
\item In cases (2), (3), and (5) we have 
$$\Phi(s)=\dfrac{C(K_2/k)}{(s-1)}+O(1),$$
with 
$$C(K_2/k)=c_1c_2c_4(\Res_{s=1}\z_{K_2'}(s))\;,$$
and
$$M(K_2/k,X)=C(K_2/k)X+O(X^{\alpha})\text{\quad for some $\alpha<1$\;.}$$
\end{itemize}
\end{corollary}

\Proof It is easy to see that when $\chi$ is not the trivial character, the
functions $F(\b,\chi,s)$ are holomorphic for $\Re(s)>1/2$, so do not 
occur in the polar part at $s=1$. On the other hand, since 
$\gothr^e(\b)\mid\gd_3$, for $\chi=1$ we have 
$F(\b,1,s)=2^{\om(\gothr^e(\b))}P(s)$, where
$P(s)=\prod_{p\in\mathcal{D}}\left(1+\dfrac{2}{\N(p)^{s}}\right)$,
so we just need to develop $P(s)$ to get the formula for the polar
part of $\Phi(s)$.

Finally, since our Dirichlet series have nonnegative and
polynomially bounded coefficients, the asymptotic results follow from a
general (and in this case easy) Tauberian theorem. For the error term
$O(X^{\alpha})$ with an explicit $\alpha<1$, we refer to the following
proposition and corollary.\fp

\begin{proposition}\label{complex}
Let $F(s)=\sum_{n=1}^\infty a_nn^{-s}$ be a Dirichlet series which is
absolutely convergent for $\Re(s)>1$, which can be extended meromorphically
to $\Re(s)>1/2$ with a pole of order $k\geq 1$ at $s=1$ and no other pole in the
strip $\frac{1}{2}<\Re(s)<1$. In addition, assume the following:
\begin{enumerate}
\item The coefficients $a_n$ are nonnegative, and for all $\varepsilon>0$
we have $$a_n\ll_{\varepsilon} n^\varepsilon\;.$$

\item $F(s)$ is a function of finite order in the vertical strip
$\frac{1}{2} < \sigma\leq 1$ : we have 
$$\vert F(\sigma+it)\vert \ll_{\varepsilon} \vert t\vert ^{\mu(\sigma) +
\varepsilon},\quad \text{when $\vert t\vert \geq 1$, for all $\varepsilon >
0$},$$where
$\mu(1)=0$, and $\mu(\sigma)$ is convex and
decreasing in the strip.

\item The integral $$\int_{0}^{1}\vert F(\sigma + it)\vert \,dt $$ is bounded
independently of $\frac{1}{2} < \sigma < \frac{1}{2} + \delta$, for some
$\delta > 0$.

\end{enumerate}
Then for all $\varepsilon>0$, we have
$$\sum_{n\leq x} a_n=\Res_{s=1}\left(F(s)\frac{x^s}{s}\right)+O(x^{\alpha+\varepsilon})\;,$$
where

\begin{equation}\label{alpha2}\alpha=1-\frac{1}{2\left(1+\mu(1/2)\right)}.\end{equation}

\end{proposition}
\Proof
Apply Perron's formula, Cauchy's residue formula and use (1) and (2) to
bound the error term.\fp

\begin{corollary}
The error term in \ref{resgen1} is  $O(x^\alpha)$, where $\alpha$ is given by
(\ref{alpha2}).
\end{corollary}
\Proof
We only need to prove that $\Phi(s)$ satisfies the hypothesis of
Proposition \ref{complex}.
For (1) we can simply refer to \cite[Lemma 6.1]{Dat-Wri} or look at the 
form of $F(\mathfrak{b},\chi,s)$, and for (2) we apply the
Phragm\'en-Lindel\"of principle.\fp

{\bf Remark.} In the case $k=\mathbb{Q}$ it is easy to show that $\mu(1/2)\leq
1/2$, so we obtain an error term $O(X^{2/3+\varepsilon})$. The previous
bound on $\mu(1/2)$ is obtained by using only the convexity bound on the
Riemann zeta function, but if we use subconvexity bounds we would
get better results.

On the other hand, if we assume the Lindel\"of hypothesis (which is for
example implied by GRH), we obtain $\mu(1/2)=0$, giving an error term
$O(X^{1/2+\varepsilon})$.

\section{Special Cases: $k=\Q$, Cases (2), (4), and (5)}\label{sec7}

We consider the case $k=\Q$, and since $\rho\notin k$ only cases (2),
(4), and (5) occur.

\subsection{Case (2): Cyclic Cubic Extensions}

\begin{proposition}\label{propcyc} We have
$$\sum_{K/\Q\text{ cyclic cubic}}\dfrac{1}{f(K)^s}=-\dfrac{1}{2}+\dfrac{1}{2}\left(1+\dfrac{2}{3^{2s}}\right)\prod_{p\equiv1\pmod3}\left(1+\dfrac{2}{p^s}\right)\;.$$
\end{proposition}

\begin{corollary} If, as above, $M(\Q,X)$ denotes the number of cyclic
cubic fields $K$ up to isomorphism with $f(K)\le X$, we have
$$\begin{array}{c}M(\Q,X)=C(\Q)X+O(X^{2/3+\varepsilon})\text{\quad where}\\
C(\Q)=\dfrac{11\sqrt{3}}{36\pi}\prod_{p\equiv1\pmod3}\left(1-\dfrac{2}{p(p+1)}\right)=0.1585282583961420602835078203575\dots\end{array}$$
\end{corollary}

\subsection{Case (4): Pure Cubic Fields}

In case (4), we have $K_2=\Q(\rho)=\Q(\sqrt{-3})$, so that $L=K_2$, and
$K/\Q$ is a pure cubic field, in other words $K=\Q(\root3\of{m})$.

\begin{proposition} We have
\begin{align*}\sum_{K/\Q\text{ pure cubic}}\dfrac{1}{f(K)^s}&=-\dfrac{1}{2}+\dfrac{1}{6}\left(1+\dfrac{2}{3^s}+\dfrac{6}{3^{2s}}\right)\prod_{p\ne3}\left(1+\dfrac{2}{p^s}\right)\\
&\phantom{=}+\dfrac{1}{3}\prod_{p\equiv\pm1\pmod{9}}\left(1+\dfrac{2}{p^s}\right)\prod_{p\not\equiv\pm1\pmod{9}}\left(1-\dfrac{1}{p^s}\right)\;,\end{align*}
where $p\not\equiv\pm1\pmod{9}$ includes $p=3$.
\end{proposition}

\begin{corollary} If, as above, $M(\Q(\sqrt{-3}),X)$ denotes the number of 
pure cubic fields $K$ up to isomorphism with $f(K)\le X$, we have
$$M(\Q(\sqrt{-3}),X)=C(\Q(\sqrt{-3}))X(\log(X)+D(\Q(\sqrt{-3}))-1)+O(X^{2/3+\varepsilon})\;,$$
where
\begin{align*}
C(\Q(\sqrt{-3}))&=\dfrac{7}{30}\prod_{p}\left(1-\dfrac{3}{p^2}+\dfrac{2}{p^3}\right)\\
&=0.066907733301378371291841632984295637501344\dots\\
D(\Q(\sqrt{-3}))&=2\ga-\dfrac{16}{35}\log(3)+6\sum_p\dfrac{\log(p)}{p^2+p-2}\\
&=3.45022279783059196279071191967111041826885\dots\;,\end{align*}
where $\ga$ is Euler's constant and the sum is over all primes including
$p=3$.
\end{corollary}

To check the validity of these constants, we note that for instance
for $X=10^{18}$ we have
\begin{align*}M(\Q(\sqrt{-3}),X)&=2937032340990444425\;,\text{\quad while}\\
C(\Q(\sqrt{-3}))X(\log(X)+D(\Q(\sqrt{-3}))-1)&=2937032340990158620\dots
\end{align*}
As already mentioned, the error is of the order of $O(X^{1/4+\eps})$ (in this
case for instance $0.22\cdot X^{1/4}\log(X)$), much smaller than 
$O(X^{2/3+\varepsilon})$  proved above, and even better than the error
term $O(X^{1/2+\varepsilon})$ that we can prove under the Lindel\"of conjecture.

\subsection{Case (5): $K_2=\Q(\sqrt{D})$ with $D\ne-3$}

In case (5), we have $K_2=\Q(\sqrt{D})$ with $D\ne-3$, so 
$L=\Q(\sqrt{D},\sqrt{-3})$. 

\begin{proposition}\label{propasym} Let $D$ be a fundamental discriminant with
$D\ne-3$, let $K_2=\Q(\sqrt{D})$, and let $r_2(D)=1$ for $D<0$ and $r_2(D)=0$
for $D>0$. There exists a function $\phi_D(s)$ holomorphic for $\Re(s)>1/2$
such that
$$\sum_{K\in\mathcal{F}(K_2)}\dfrac{1}{f(K)^s}=\phi_D(s)+\dfrac{3^{r_2(D)}}{6}L_3(s)\prod_{\lgs{-3D}{p}=1}\left(1+\dfrac{2}{p^s}\right)\;,\quad\textrm{
where}$$

$$L_3(s)=\begin{cases} 
1+2/3^{2s}&\text{\quad if $3\nmid D$,}\\
1+2/3^s&\text{\quad if $D\equiv3\pmod9$,}\\
1+2/3^s+6/3^{2s}&\text{\quad if $D\equiv6\pmod9$.}
\end{cases}$$
\end{proposition}

\Proof If we denote by $\phi_D(s)$ the contribution of the nontrivial 
characters in Theorem \ref{mainth1} it is clear that $\phi_D(s)$ is a
holomorphic function for $\Re(s)>1/2$, so it is sufficient to consider the
contribution of the trivial characters $\Phi_0(s)$. We consider the three cases 
separately and, with similar notations and computations as in the examples above, we
get:
$$\begin{array}{ll}
 \Phi_0(s)=\dfrac{3^{r_2(D)}}{6}\left(1+\dfrac{2}{3^{2s}}\right)\prod_{\lgs{-3D}{p}=1}\left(1+\dfrac{2}{p^s}\right)
 & \textrm{if }3\nmid D\\
\Phi_0(s)=\dfrac{3^{r_2(D)}}{6}\left(1+\dfrac{2}{3^s}\right)\prod_{\lgs{-3D}{p}=1}\left(1+\dfrac{2}{p^s}\right)
& \textrm{if }D\equiv3\pmod9\\
\Phi_0(s)=\dfrac{3^{r_2(D)}}{6}\left(1+\dfrac{2}{3^s}+\dfrac{6}{3^{2s}}\right)\prod_{\lgs{-3D}{p}=1}\left(1+\dfrac{2}{p^s}\right)
& \textrm{if }D\equiv6\pmod9,D\ne-3\;.\end{array}$$
\fp

\smallskip
\begin{corollary}\label{corasym} Set $D'=-3D$ if $3\nmid D$ and $D'=-D/3$ if
$3\mid D$, and denote as usual by $\chi_{D'}$ the character $\lgs{D'}{.}$. Then
if $D\ne-3$ is a fundamental discriminant we have
\begin{align*}M(\Q(\sqrt{D}),X)&=C(\Q(\sqrt{D}))X+O(X^{2/3+\varepsilon})\text{ where}\\
C(\Q(\sqrt{D}))&=\dfrac{3^{r_2(D)}\ell_3L(\chi_{D'},1)}{\pi^2}\prod_{p\mid D'}\left(1-\dfrac{1}{p+1}\right)\prod_{\lgs{D'}{p}=1}\left(1-\dfrac{2}{p(p+1)}\right)\;,\end{align*}
where
$$\ell_3=\begin{cases}
11/9&\text{\quad if $3\nmid D$,}\\
5/3&\text{\quad if $D\equiv3\pmod9$,}\\
7/5&\text{\quad if $D\equiv6\pmod9$.}
\end{cases}$$
\end{corollary}

Note that $L(\chi_{D'},1)$ is given by Dirichlet's class number formula, in
other words with standard notation, $L(\chi_{D'},1)=2\pi h(D')/(w(D')\sqrt{|D'|})$ if $D'<0$ and $L(\chi_{D'},1)=2h(D')R(D')/\sqrt{D'}$ if $D'>0$.

\bigskip

\myproof{ of Theorem \ref{thresq} (2)} We now show how to modify the above
formulas so as to obtain the formula given in the theorem. By Propositions
\ref{propcyc} and \ref{propasym} we can write
$$\Phi_D(s)=\phi_D(s)+g(K'_2)\dfrac{3^{r_2(D)}}{6}L_3(s)\prod_{p\ne3}\left(1+\dfrac{a_{K'_2}(p)}{p^s}\right)\;,$$
where $g(K'_2)=1$ unless $D=1$, in other words $K'_2=\Q(\sqrt{-3})$, in which
case $g(K'_2)=3$. Thus,
$$\dfrac{\Phi_D(s)}{(1-1/3^s)\z(s)}=\psi_D(s)+g(K'_2)\dfrac{3^{r_2(D)}}{6}L_3(s)\prod_{p\ne3}\left(1+\dfrac{a_{K'_2}(p)}{p^s}\right)\left(1-\dfrac{1}{p^s}\right)\;,$$
where $\psi_D(s)=\phi_D(s)/((1-1/3^s)\z(s))$.
When $s$ tends to $1$, $\psi_D(s)$ tends to $0$, the left-hand side tends to
a limit, and it is easy to see that the right-hand side tends to a 
semi-convergent Euler product. Thus, if we set 
$P(K'_2)=\prod_{p\ne 3}((1+a_{K'_2}(p)/p)(1-1/p))$, we have
$$C(\Q(\sqrt{D})=\Res_{s=1}\Phi_D(s)=g(K'_2)\dfrac{1}{3^{2-r_2(D)}}L_3(1)P(K'_2)
=g(K'_2)\dfrac{c_3(K'_2)}{3^{3+r_2(K'_2)}}P(K'_2)\;,$$
where $c_3(K'_2)$ is given in the theorem, since the different cases for
$L_3(1)$ correspond to the different splittings of $3$ in $K'_2/\Q$.\fp

\subsection{An Exact Result when $D<0$ and $3\nmid h(D)$}\label{schol}

It is interesting to note that when $D<0$ and $3\nmid h(D)$, one can prove that
nontrivial characters do not occur in the above formulas, so that 
$\phi_D(s)=0$, thus giving exact formulas for the Dirichlet series.

\begin{proposition} Assume that $K_2=\Q(\sqrt{D})$ with $D<0$, $D\ne-3$, and
$3\nmid h(D)=|Cl(K_2)|$. Then for any ideal $\b\in\mathcal{B}$ occurring in
the sum of Theorem \ref{mainth1}, the group $G_{\b}=(Cl_{\b}(L)/Cl_{\b}(L)^3)[T]$ is trivial.\end{proposition}

\Proof An important theorem of Scholz (\cite{Sch}) says that if $D<0$ is
a negative fundamental discriminant different from $-3$ we have
$$0\le \rk_3(Cl(\Q(\sqrt{D})))-\rk_3(Cl(\Q(\sqrt{-3D})))\le 1$$
and that $ \rk_3(Cl(\Q(\sqrt{D})))=\rk_3(Cl(\Q(\sqrt{-3D})))$ if and
only if the fundamental unit $\eps$ of $\mathbb{Q}(\sqrt{3D})$ is not $3$-primary, in other words if and only if $\eps$ is not
a cube modulo $3\sqrt{-3}\Z_L$, where $L=\Q(\sqrt{D},\sqrt{-3})$. Since in our case we assume
that $\rk_3(Cl(\Q(\sqrt{D})))=0$, it follows that we also have
$\rk_3(Cl(\Q(\sqrt{-3D})))=0$ and that $\eps$ is not a cube modulo 
$3\sqrt{-3}\Z_L$.

We now consider the exact sequence of $\F_3[G]$-modules already used
above in the computation of $f_{\al_0}(\b)$:
$$1\LR S_{\b}(L)[T]\LR S_3(L)[T]\LR \dfrac{Z_{\b}}{Z_{\b}^3}[T]\LR\dfrac{ Cl_{\b}(L)}{Cl_{\b}(L)^3}[T]\LR\dfrac{Cl(L)}{Cl(L)^3}[T]\LR 1\;.$$
By Hasse's formula giving the class number of biquadratic number fields
(\cite{Has}), we have $|Cl(L)|=2^{-j}|Cl(K_2)||Cl(K_2')|$ with $j=0$ or $1$,
so in particular by Scholz's theorem we deduce that $3\nmid|Cl(L)|$.
We thus have the exact sequence
$$1\LR S_{\b}(L)[T]\LR S_3(L)[T]\LR \dfrac{Z_{\b}}{Z_{\b}^3}[T]\LR G_{\b}\LR 1\;.$$
In addition, also since $3\nmid|Cl(L)|$, $S_3(L)$ is an $\F_3$-vector space
of dimension $r_1(L)+r_2(L)=2$, generated by the classes modulo cubes of
$\rho$ and a fundamental unit $\eps$ of $K'_2=\Q(\sqrt{-3D})$. The action of 
$\tau$ and $\tau_2$ is given by $\tau(\rho)=\rho^{-1}$, $\tau_2(\rho)=\rho$,
$\tau(\eps)=\pm\eps^{-1}$, $\tau_2(\eps)=\pm\eps^{-1}$ (where
$\pm=\N_{K'_2/\Q}(\eps)$), and modulo cubes the $\pm$ signs disappear.
Since $T=\{\tau+1,\tau_2+1\}$, it follows that $S_3(L)[T]$ is a $1$-dimensional
$\F_3$-vector space generated by the class of $\eps$. 

Since $G_{\b}$ maps surjectively onto $G_{\b'}$ for $\b'\mid\b$, it is 
sufficient to consider $\b=3\sqrt{-3}$. In that case, we have seen that
$|(Z_{\b}/Z_{\b}^3)[T]|=3$ in all cases, and since we have just shown that
$|S_3(L)[T]|=3$, by the above exact sequence it follows that $G_{\b}$ is 
trivial if and only if $S_{\b}(L)[T]$ is trivial, hence by definition if
and only if $\eps$ is not congruent to a cube modulo $\b_z=3\sqrt{-3}\Z_L$,
which is exactly the second statement of Scholz's theorem, proving the
proposition.\fp

\smallskip

\noindent
{\bf Remark.}
The same proof shows the following result for $D>0$:
if $D>0$ and $3\nmid h(D')$, where as usual $D'=-3D$ if $3\nmid D$ and 
$D'=-D/3$ if $3\mid D$, then $G_{\b}$ is canonically isomorphic to
$(Z_{\b}/Z_{\b}^3)[T]$, hence has order $1$ unless $\b=3\sqrt{-3}$ or
$3\nmid D$ and $\b=3\Z_L$, in which case it has order $3$.

\medskip

\begin{corollary} Under the same assumptions, we have the following simple 
result:
$$\sum_{K\in\mathcal{F}(K_2)}\dfrac{1}{f(K)^s}=-\dfrac{1}{2}+\dfrac{1}{2}L_3(s)\prod_{\lgs{-3D}{p}=1}\left(1+\dfrac{2}{p^s}\right)\;,$$
where
$$L_3(s)=\begin{cases}
1+2/3^{2s}&\text{\quad if $3\nmid D$,}\\
1+2/3^s&\text{\quad if $D\equiv3\pmod9$,}\\
1+2/3^s+6/3^{2s}&\text{\quad if $D\equiv6\pmod9$.}
\end{cases}$$
\end{corollary}

\enddocument